\documentclass[11pt,a4paper]{article}

\usepackage{amssymb,latexsym,amsmath,epsfig,amsthm} 
\usepackage{url}
\usepackage{lmodern}
\usepackage{caption}
\usepackage{adigraph}
\usepackage{graphicx}
\usepackage{subcaption}
\usepackage{diagbox}
\usepackage{amscd}
\usepackage{verbatim}
\usepackage{tikz}
\usepackage{youngtab}
\usepackage{ytableau}
\usepackage{here}
\usepackage{cite}
\usepackage{hyperref}

\newcommand{\m}{\boldsymbol m}
\newcommand{\x}{\boldsymbol x}
\newcommand{\y}{\boldsymbol y}

\newcommand\Zz{{\mathbb Z_{\geq 0}}}

\newcommand\Mex{{\mathrm{mex}}}

\newcommand\Pp{${\mathcal{P}}$-position}
\newcommand\Pps{${\mathcal{P}}$-positions}
\newcommand\Np{${\mathcal{N}}$-position}

\newcommand\RMi{\,{\rm i}\,}
\newcommand\RMii{\,{\rm ii}\,}
\newcommand\RMiii{\,{\rm iii}\,}

\makeatother

\theoremstyle{plain}
\newtheorem{theorem}{Theorem}[section]
\newtheorem{remark}[theorem]{Remark}
\newtheorem{proposition}[theorem]{Proposition}
\newtheorem{lemma}[theorem]{Lemma}

\theoremstyle{definition}
\newtheorem{definition}[theorem]{Definition}
\numberwithin{equation}{section}
\numberwithin{table}{section}
\numberwithin{figure}{section}

\newcommand{\ifEmpty}[3]{\setbox0=\hbox{#1\unskip}\ifdim\wd0=0pt{#2}\else{#3}\fi}
\newcommand{\nothingIfEmpty}[2]{\ifEmpty{#1}{}{#2}}
\newcommand{\showIfNotEmpty}[1]{\nothingIfEmpty{#1}{#1}}

\newcommand{\TalkTitle}{Combinatorial Games and the Golden Ratio on Digraphs}
\newcommand{\TalkPresenter}{
Tomoaki Abuku\footnote{The first author is supported by JSPS Kakenhi 22K13953.},
Hiroki Inazu, 
Hiyu Inoue,
Shun-ichi Kimura\footnote{The fourth author was supported by JSPS Kakenhi 23K03071.},
Koki Suetsugu\footnote{The fifth author is supported by The INOUE ENRYO Memorial Grant, TOYO University.},
Kosaku Watanabe,
Takahiro Yamashita\footnote{The seventh author was supported by JST SPRING, Grant Number JPMJSP2132.}
}
\newcommand{\TalkEmail}{
buku3416@gmail.com,
d220585@hiroshima-u.ac.jp,
hiyuuinoue@gmail.com,
skimura@hiroshima-u.ac.jp,
suetsugu.koki@gmail.com,
d206934@hiroshima-u.ac.jp,
d236676@hiroshima-u.ac.jp
}

\begin{document}
\begin{center}
 {\Large\textsc{\TalkTitle}}

 \textbf{\TalkPresenter}

 \showIfNotEmpty{\texttt{\TalkEmail}}



 \bigskip\bigskip
\end{center}

\begin{abstract}
We introduce a new combinatorial game called Triangle Game. 
In this game, a directed $3$-cycle graph is given, and stones are placed on each vertex. 
The player chooses a directed edge and takes at least one stone from the initial vertex. At the same time, the player is allowed to return some stones to the terminal vertex of the edge, as long as the total number of stones decreases. 
We describe the set of \Pps~under both normal play and mis\`ere play.
The golden ratio $\phi=\dfrac{1+\sqrt{5}}{2}$ plays an essential role in our description. We also show that Triangle Game is tame.
\end{abstract}


\section{Introduction}

In this paper, we study combinatorial games, namely, two-player games involving neither chance nor hidden information. One of the most fascinating aspects of combinatorial game theory is that the analysis of seemingly simple rules often reveals unexpected mathematical structures and constants.

A classical example is Wythoff Nim. A position of this game is represented by a pair of non-negative integers $(x,y)$. On each turn, a player may decrease either coordinate by an arbitrary positive integer or decrease both coordinates by the same positive integer. The player who moves to $(0,0)$ wins. Remarkably, the golden ratio arises in the complete characterization of the $\mathcal{P}$-positions (positions in which the previous player has a winning strategy). More precisely, a position is a $\mathcal{P}$-position if and only if it is of the form
\[
(\lfloor n\phi\rfloor,\lfloor n\phi\rfloor+n)
\quad\text{or}\quad
(\lfloor n\phi\rfloor+n,\lfloor n\phi\rfloor)
\]
for some non-negative integer $n$, where
\[
\phi=\frac{1+\sqrt{5}}{2}
\]
is the golden ratio~\cite{Wythoff}.

The golden ratio also appears in Euclid Nim, a two-heap impartial game introduced in~\cite{CD69}. A position is represented by a pair $(a,b)$ of positive integers. If $a>b$, a player may remove $kb$ stones from the larger heap, where $k$ is a positive integer such that $a-kb\geq 0$; the case $b>a$ is defined symmetrically. A position $(a,b)$ with $a>b>0$ is a $\mathcal{P}$-position if and only if
\[
a<\phi b.
\]

The appearance of the golden ratio in game analysis is quite rare, and only a small number of examples, including those described above, are known.  

The main contribution of this paper is to provide a new natural example of an impartial combinatorial game in which the golden ratio appears in the structure of the \(\mathcal P\)-positions.
More
precisely, in Triangle Game, which will be defined in the next section, the golden ratio does not appear in the
rules themselves; rather, it arises in the explicit description of the
\(\mathcal P\)-positions under both normal play (last player wins) and mis\`ere play (last player loses).

We also compare the normal and mis\`ere Sprague--Grundy values. Although we were unable to obtain a closed-form expression that determines the Sprague--Grundy value of every position, we prove that the ruleset is tame.

The remainder of this paper is organized as follows. In the remainder
of this section, we introduce the basic terminology and results
concerning impartial games, Sprague--Grundy values, and tame rulesets.
In Section~\ref{sec:Digraph_Nim}, we define Triangle Game and prove
our main results concerning its \(\mathcal P\)-positions and
tameness. Finally, in Section~\ref{Sec:final}, we present several open
questions and directions for future research.

\subsection{Impartial Games}
\label{sec:def_impartial}

\begin{definition}[Impartial game]
\label{def:impartial}
An \textbf{impartial game} is a triple $\Gamma=(M, f, w)$, where $M$ is the set of game positions, $f: M\to Pow(M)$ is the option map, with $Pow(M)$ the set of subsets of $M$, for $\m \in M$, $\m'\in f(\m)$ means that the player can make the move from $\m$ to $\m'$. We also denote this situation simply by $\m\to \m'$, following the usual notation. The symbol $w$ is the rule to determine the winner, we treat only two cases, say $w\in \{\textbf{Normal}, \textbf{Mis\`ere}\}$.
When $w=\textbf{Normal}$, we play the game under normal play convention, which means that the last player wins, and when $w=\textbf{Mis\`ere}$, we play the game under mis\`ere play convention, which means that the last player loses.
When the play convention is irrelevant, we simply write $\Gamma=(M, f)$.
We also assume that the game terminates in finitely many steps. To be more precise, starting from each $\m \in M$, there exists a non-negative integer $\ell(\m)$ such that the length of the plays $\m=m_0 \to m_1 \to \cdots \to m_k$ is bounded as $k\le \ell(\m)$. We say that $\m\in M$ is a \textbf{terminal position} when $f(\m)=\emptyset$, or equivalently when we can take $\ell(\m)=0$.
We write the set of all terminal positions as $\mathcal{E}$.
\end{definition}

\begin{definition}[\Np~and \Pp]\label{def:PpNp}
Let $\Gamma=(M,f,w)$ be an impartial game.
We call a position $\m\in M$ an \Np~if the current player has a winning strategy.
We call a position $\m\in M$ a \Pp~if the previous player has a winning strategy.
\end{definition}

\begin{proposition}
\label{prop:PpNp}
  Let $\Gamma=(M, f, w)$ be an impartial game. For all $\m\in M$, we can determine who has a winning strategy as follows:
  \begin{enumerate}
  \item[(i) ] If $\m$ is a terminal position, when $w=\textbf{Normal}$, $\m$ is a \Pp, and if $w=\textbf{Mis\`ere}$, $\m$ is an \Np.
  \item[(ii) ] 
  If there exists $\m' \in f(\m)$ such that $\m'$ is a \Pp, then $\m$ is an \Np.
  \item[(iii) ] 
  If $f(\m)\neq\emptyset$ and for any $\m'\in f(\m)$, $\m'$ is an \Np, then $\m$ is a \Pp.
  \end{enumerate}
  In particular, any position $\m\in M$ is either a \Pp~or an \Np.
\end{proposition}
This standard characterization follows by induction on the game tree; see, for example, Siegel~\cite{S}.

\begin{remark}
Note that Proposition \ref{prop:PpNp} states that every non-terminal \Np~has an option that is a \Pp, and the winning strategy is to move to such a \Pp.
\end{remark}

Furthermore, a finer classification is possible using a parameter called the Sprague--Grundy value.
Let $\Zz$ be the set of non-negative integers.

\begin{definition}[Mex]
\label{mex}
Let $X$ be a finite subset of non-negative integers.
Then the \textbf{minimal excluded number $\mathrm{mex}(X)$} is 
$$\mathrm{mex}(X)=\min(\mathbb{Z}_{\geq 0} \setminus X).$$
\end{definition}

\begin{definition}[Sprague--Grundy value under normal play]
\label{nim_value_normal}
Let $\Gamma = (M,f)$ be an impartial game.
For any $\m \in M$, the \textbf{Sprague--Grundy value} for normal play $g^+(\m)$ is
$$g^+(\m)= \mathrm{mex}(\{g^+({\m}') \mid \m'\in f(\m)\}).$$
\end{definition}

\begin{remark}
    Note that $g^+(\m) = 0~\text{if}~\m\in \mathcal{E}~\text{since}~\{g^+(\m')\mid \m' \in f(\m)\}$ is an empty set.
\end{remark}

\begin{proposition}
\label{Sprague--Grundy}
Let $\Gamma = (M, f,\textbf{Normal})$ be an impartial game under normal play.
For any $\m \in M$, we have the following properties.
 \begin{enumerate}
    \item[\text{(1)} ] $\m$ is a \Pp~if and only if $g^+(\m) = 0$;
    \item[\text{(2)} ] $\m$ is an \Np~if and only if $g^+(\m) \neq 0$.
 \end{enumerate}
\end{proposition}

\begin{definition}[Disjunctive sum]
\label{def:disjunctive_sum}
Let $X=(M_1,f_1)$ and $Y=(M_2,f_2)$ be impartial games.
The \textbf{disjunctive sum} $X + Y$ is an impartial game $(M,f)$ where the set of game positions is $M=M_1 \times M_2$ and for any $\x\in M_1$ and $\y \in M_2$, the option map $f(\x+\y)$ is
\begin{eqnarray*}
       f(\x+\y) = \{(\x'+ \y)\in M\, |\, \x' \in f_1(\x)\}  \cup \{(\x+ \y')\in M\, |\, \y' \in f_2(\y)\}.
\end{eqnarray*}
\end{definition}

\begin{remark}
     Definition \ref{def:disjunctive_sum} means on each turn, a player chooses one of the two component game positions and makes a legal move from that position, leaving the other unchanged.
\end{remark}

\begin{theorem}[Sprague--Grundy Theorem \cite{SP35,GR39}]
\label{thm:Sprague--Grundygeneral}
Let $X + Y=(M,f, \textbf{Normal})$ be an impartial game under normal play where $X=(M_1,f_1),Y=(M_2,f_2)$. For any $\x\in M_1$ and for any $\y \in M_2$, we have
$$g^+(\x+\y)= g^+(\x) \oplus  g^+(\y),$$
where $\oplus$ denotes the bitwise XOR operation on binary representations.
\end{theorem}

\begin{definition}[Sprague--Grundy value under mis\`ere play]
\label{nim_value_misere}
Let $\Gamma = (M,f)$ be an impartial game.
For any $\m \in M$, the \textbf{Sprague--Grundy value} for mis\`ere play $g^-(\m)$ is
$$g^-(\m)=\left\{
\begin{array}{cc} 
1& (\text{if}~\m\in \mathcal{E})\\ 
\mathrm{mex}(\{g^-({\m}') \mid \m'\in f(\m)\}) & (\text{if}~\m\notin \mathcal{E}).
\end{array}
\right.
$$
\end{definition}
\begin{proposition}
\label{Sprague--Grundy_misere}
Let $\Gamma = (M, f, \textbf{Mis\`ere})$ be an impartial game under mis\`ere play.
For any $\m \in M$, we have the following properties.
 \begin{enumerate}
    \item[\text{(1)} ] $\m$ is a \Pp~if and only if $g^-(\m) = 0$;
    \item[\text{(2)} ] $\m$ is an \Np~if and only if $g^-(\m) \neq 0$.
 \end{enumerate}
\end{proposition}

\begin{definition}[Tame\cite{C76}]
\label{def_tame}
A position \(\m\) is called an \textbf{$(i,j)$-position} if
\[
g^+(\m)=i
\quad\text{and}\quad
g^-(\m)=j.
\]
An impartial ruleset is called \textbf{tame} if every position is
a \((0,1)\)-position, a \((1,0)\)-position, or a
\((k,k)\)-position for some \(k\in\mathbb Z_{\ge0}\).
\end{definition}

In general, the Sprague--Grundy values under mis\`{e}re play behave quite differently from those under normal play. However, once a game is known to be tame, the two are closely related. For this reason, determining whether a game is tame has attracted considerable interest
(for the details of tame games, see \cite{C76,S}).

\section{Main results}\label{sec:Digraph_Nim}

We now determine the \Pps~of Triangle Game under both normal play and mis\`ere play. We also show that Triangle Game is tame.

\subsection{The rules of Triangle Game}
\begin{definition}[Triangle Game]
\label{def:TriangleGame}

\textbf{Triangle Game} is an impartial game played on a directed graph $G = (V, E)$ with the set of vertices $V=\{X, Y, Z\}$ and the set of directed edges $E=\{(X, Y), (Y, Z), (Z, X)\}$. Stones are placed on each vertex.
On each turn, the player chooses an edge $(V_s, V_t) \in E$, and takes $i\ge1$ stones from $V_s$, and may return $j$ of them to $V_t$, where $0\le j<i$. Thus, the total number of stones strictly decreases.

To be more concrete, this is an impartial game $\Gamma=(M, f, w)$ with 
\begin{eqnarray*}
  M &=& (\Zz)^3 =\{(x, y, z) \mid x, y, z\in \Zz\},\\
  f((x,y,z)) &=&\{(x-i,y+j,z) \mid 1 \leq i \leq x, 0 \leq j <i \} \\
           & &\cup~ \{(x,y-i,z+j) \mid 1 \leq i \leq y, 0 \leq j <i \} \\
           & &\cup ~\{(x+j,y,z-i) \mid 1 \leq i \leq z, 0 \leq j <i \},
  \end{eqnarray*}
\end{definition}
where $x, y$ and $z$ are the numbers of stones on the vertices $X, Y$ and $Z$ respectively.

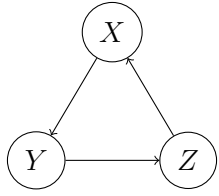
\begin{figure}[htb]
\begin{center}
\begin{tikzpicture}[auto]
\node[shape=circle, draw] (a) at (1, 0) {$X$};
\node[shape=circle, draw] (b) at (0, -1.7) {$Y$};
\node[shape=circle, draw] (c) at (2, -1.7) {$Z$};
\draw[->,very thick] (b) to (c);
\draw[->,very thick] (a) to (b);
\draw[->,very thick] (c) to (a);
\end{tikzpicture}
\caption{\textbf{Graph of Triangle Game}}
\label{fig:graphG}
\end{center}
\end{figure}

Figure \ref{fig:graphG} represents the graph of Triangle Game.

\begin{remark}
    For the games treated in this paper, we may take $\ell(\m)$ to be the total number of stones in the position $\m\in M$.
\end{remark}

\begin{remark}
    The set of all terminal positions for Triangle Game is
    $$\mathcal{E}=\{(0,0,0)\}.$$
\end{remark}

\subsection{Triangle Game under normal play}

\begin{theorem}
\label{Thm:TrianglePps_normal}
Let $S$ be the set of $\mathcal{P}$-positions of Triangle Game under normal play.
Then $S$ is given by
\begin{eqnarray*}
\label{Eqn:f_functionpart1Pps}
S=\{
(a,b,c),(b,c,a),(c,a,b)\mid a,b,c\in\Zz,\ a=b+c,\ b\ge \phi c\},
\end{eqnarray*}
where $\phi=\frac{1+\sqrt{5}}{2}$.
\end{theorem}

\begin{remark}
Since \(\phi\) is irrational, the condition \(b\ge \phi c\) is equivalent to \(b>\phi c\) whenever \(c>0\). 
As $(0, 0, 0)$ is a \Pp, we need the equality here.
\end{remark}

\begin{lemma}
\label{Lemma_Hi}
Let $x, y, z$ be positive integers satisfying $x=y+z$ with $x>y\geq z>0$.
Then, we have the following properties:
\begin{enumerate}
    \item [(\RMi)] $x>\phi y \iff y<\phi z$;
    \item [(\RMii)] $x<\phi y \iff y>\phi z$.
\end{enumerate}
\end{lemma}
\begin{proof}
We only show property (\RMi), as (\RMii) is equivalent to (\RMi) because $x,y$ and $z$ are positive integers and $\phi$ is irrational.
\begin{eqnarray*}
             x>\phi y &\iff& \dfrac{x}{y}>\phi \\
                      &\iff& \dfrac{y+z}{y}>\phi \\
                      &\iff& y+z>\phi y \\
                      &\iff& z>(\phi-1)y \\
                      &\iff& z>\dfrac{1}{\phi}y~\left(\text{because}~\phi = 1+ \dfrac{1}{\phi}\right) \\
                      &\iff& \phi z>y 
\end{eqnarray*}
\end{proof}

\begin{proof}[Proof of Theorem \ref{Thm:TrianglePps_normal}]
By Proposition~\ref{prop:PpNp}, we need to show that
\begin{enumerate}
\item[(\RMi)] $\mathcal{E} \subset S$;
\item[(\RMii)] For any $(x,y,z)\in S$ and its option $(x',y',z')\in f((x,y,z))$, we have $(x',y',z')\notin S$;
 \item[(\RMiii)] For any  $(x,y,z)\notin S$, there is an option $(x',y',z')\in f((x,y,z))$ such that $(x',y',z')\in S$.
\end{enumerate}

It is obvious that (\RMi) holds.

For (\RMii), we take $(x,y,z)\in S$. 
Without loss of generality, we assume $x=y+z$ with $y \geq \phi z$. 
Let $(x',y',z')$ be an option for $(x,y,z)$.
If we take stones from $y$ or $z$, we have $y'+z'<y+z=x\leq x'$.
Then, we have $x'>y'+z'$ and hence $(x',y',z') \notin S$. If we take stones from $x$, then we have $x>x'$ and $y'\geq y\geq z=z'$.
As $y'\ge z'$, the value of $x'$ can be one of the following cases below:
\begin{enumerate}
    \item[(1)] $x'\ge y'\ge z'$; 
    \item[(2)] $y'>x'>z'$;
    \item[(3-1)] $y' \ge z'\ge x'$ with  $z'>0$;
    \item[(3-2)] $y' \ge z'= x'=0$.
\end{enumerate}
We show that $(x',y',z')\notin S$ for these cases one by one.
\begin{itemize}
    \item[(1)] When $x'\geq y'\geq z'$, since we have $x' < y'+z'$, $(x',y',z') \notin S$. 
    \item[(2)] When $y'> x'> z'$, since the order does not match the assumption (when $(x', y', z')\in S$ with $y'>x', z'$, we need $z'\ge \phi x'\ge x'$), hence $(x',y',z') \notin S$.
    \item[(3-1)] 
    When $y'\ge z'\ge x'$ with $z'>0$, we have $y'\ge y>\phi z=\phi z'$. In order that $(y', z', x')\in S$, we need to have $y'=z'+x'$, but then, by Lemma \ref{Lemma_Hi}, $y'>\phi z'$ implies $z'<\phi x'$, contradicting $(y', z', x')\in S$. Hence $(x', y', z')\notin S$.
    \item[(3-2)] When $y'\ge z'= x'=0$,     if $(y', z', x')\in S$, we have $y'=z'+x'=0$.
    Because we took stones from $x$, we have $z=z'=0$ and  $y\le y'=0$ which implies $y=0$ and hence $x=y+z=0$. We cannot take stones from $(x, y, z)=(0, 0, 0)$, so this case does not happen.
\end{itemize}

For (\RMiii), we take $(x,y,z)\notin S$. Without loss of generality, we assume $x\geq y,z$. We divide into the following cases:
\begin{enumerate}
    \item[(1)] $z\geq y$;
    \item[(2)] $\phi z > y > z$;
    \item[(3)] $y \geq \phi z$ and $x>y+z$;
    \item[(4)] $y \geq \phi z$ and $x<y+z$.
\end{enumerate}
We construct $(x,y,z)\to (x',y',z')\in S$ for these cases one by one.
\begin{itemize}
    \item[(1)] In this case, we have $x\geq z\geq y$, then if $(x, y, z)$ were in $S$, then $(x, y, z)$ would be $(x, 0, x)$, that satisfy $z=x+y$ and $x\ge \phi y$.
Hence by $(x, y, z)\notin S$ we may assume either $x>z$ or $y \neq 0$, and in either case, we have $x>z-y$, so we can make the move $(x, y, z)\to (0, z, z) \in S$.
\end{itemize}
    For the rest of the proof, we may assume  $x\ge y>z$. 
\begin{itemize}
    \item[(2)] When $\phi z>y>z$, applying Lemma \ref{Lemma_Hi} to $y=z+(y-z)$, we obtain $\phi(y-z)<z$, which implies $x>y-z$ because $x>z$ and $\phi(y-z)>y-z$, we can make the move $(x, y, z)\to (y-z, y, z)\in S$.
    \item[(3)] When $y \geq \phi z$ and $x>y+z$, we can make the move $(x,y,z)\to(y+z,y,z)\in S$.
    \item[(4)] When $y \geq \phi z$ and $x<y+z$, we can make the move $(x,y,z)\to(x,0,x)\in S$ because $y>x-z$.
\end{itemize}
\end{proof}

\subsection{Triangle Game under mis\`ere play}
\begin{theorem}
\label{Thm:TrianglePps_misere}
Let $S^-$ be the set of $\mathcal{P}$-positions of Triangle Game under mis\`ere play.
Then $S^-$ is given by
$$S^-=S^-_1\cup S^-_2,$$
where
\begin{eqnarray*}
  S^-_1 &=&\{(1,0,0),(0,1,0),(0,0,1),(1,1,1)\}, \\
  S^-_2 &=& \{(a,b,c),(b,c,a),(c,a,b)
\mid a,b,c\in\mathbb Z_{\ge0},\ a=b+c\ge2,\ b\ge \phi c\}.
\end{eqnarray*}
\end{theorem}

\begin{proof}
By Proposition~\ref{prop:PpNp}, we need to show that
\begin{enumerate} 
 \item[(\RMi)]
 $\mathcal{E}\cap S^-=\varnothing$;
 \item[(\RMii)]
 For any $(x,y,z)\in S^-$ and its option $(x',y',z')\in f((x,y,z))$, we have $(x',y',z')\notin S^-$;
 \item[(\RMiii)] For any  $(x,y,z)\notin (S^-\cup\mathcal{E})$, there is an option $(x',y',z')\in f((x,y,z))$ such that $(x',y',z')\in S^-$.
\end{enumerate}

It is obvious that (\RMi) holds.

For (\RMii), we take $(x,y,z)\in S^-$. 
If $(x,y,z)\in S^-_1$, without loss of generality, we may assume to take a stone from $(1,0,0),(1,1,1)$.
When $(x,y,z)=(1,0,0)$, we can only move to $(0,0,0)\notin S^{-}$. When $(x,y,z)=(1,1,1)$, without loss of generality, we may take the stone from the vertex $X$, then we can only move to $(0,1,1)\notin S^{-}$.
If $(x, y, z)\in S^-_{2}$, without loss of generality, we may assume $x=y+z\ge 2$ with $y\ge \phi z$. Since $S^-_2\subset S$, where $S$ is the set of \Pps~under normal play, from Theorem \ref{Thm:TrianglePps_normal}, it is enough to show that when we move $(x, y, z)\to (x', y', z')$, then $(x', y', z')\notin S^-_1$. Suppose, to the contrary, that $(x', y', z')\in S^-_1$, then we take stones from $x$ as $2\le x$ and $x>1\ge x'$, hence we have $y\le y'\le 1$ and $z'=z$. As $2\le y+z\le y'+z'\le 2$, we have $y=z=1$, contradicting $y\ge \phi z=\phi>1$. Thus, $(x', y', z')\notin S^-$.

For (\RMiii), we take $(x,y,z)\notin (S^-\cup\mathcal{E})$. 
Without loss of generality, we may assume $x\geq y,z$. 
We divide into the following cases:
\begin{enumerate}
    \item[(1)] $y,z\in\{0,1\}$;
    \item[(2)] $y\geq2$ or $z\geq2$.
\end{enumerate}
We construct $(x,y,z)\to (x',y',z')\in S^-$ for these cases one by one.
\begin{enumerate}
    \item[(1)] In this case, we have the following positions and moves:
    $(x,0,0) \to (1,0,0),~(x,0,1)\to(0,0,1),~(x,1,0)\to(0,1,0),~(x,1,1)\to(1,1,1),$
    where the legality of each move follows from the assumption $(x,y,z)\notin (S^-\cup\mathcal{E})$.
    \item[(2)]
    In this case, we have $\max(\{x,y,z\})\geq2$.
    Since $(x,y,z)\notin (S^-\cup\mathcal{E})$, we have $(x,y,z)\notin S^-_2$.
    As $S_2^-=\{(a,b,c)\in S\mid \max(\{a,b,c\})\geq2\}$, it follows that $(x,y,z)\notin S$.
    From Theorem \ref{Thm:TrianglePps_normal}, we can make the move $(x,y,z)\to (x',y',z')\in S$.
If we take stones from $x$, then $y'\geq y$
and $z'= z$, so either $y' \geq 2$ or $z' \geq 2$. 
If we take stones from $y$, then $x'= x \geq 2$. 
If we take stones from $z$, then $x'= x + j \geq x \geq 2$. 
Therefore, $\max(\{x', y', z'\})\geq 2$.
Thus, $(x',y',z')\in S^-_2$.
\end{enumerate}
\end{proof}

\subsection{Triangle Game is tame}

Let
$$S_1^+=\{(0,0,0),(1,1,0),(1,0,1),(0,1,1)\},$$
and recall that

$$S_1^-=\{(1,0,0),(0,1,0),(0,0,1),(1,1,1)\}.$$
We set
$$A=S_1^+\cup S_1^-=
\{(x,y,z)\in M\mid x,y,z\in\{0,1\}\}.$$

Recall that, for Triangle Game, $\ell(\m)$ is the total number of stones in the position $\m$. For every $\m\in A$,
the number $\ell(\m)$ is even if $\m\in S_1^+$, whereas it is odd if $\m\in S_1^-$.

\begin{lemma} \label{gyakuten}
For every $\m\in A$,
\[
\bigl(g^+(\m),g^-(\m)\bigr)
=
\begin{cases}
(0,1),&\text{if }\ell(\m)\text{ is even},\\
(1,0),&\text{if }\ell(\m)\text{ is odd}.
\end{cases}
\]
In particular, every position in $S_1^+$ is a $(0,1)$-position,
and every position in $S_1^-$ is a $(1,0)$-position.
\end{lemma}

\begin{proof}
By cyclic symmetry, it is enough to consider the four positions
\[
(0,0,0),\qquad
(1,0,0),\qquad
(1,1,0),\qquad
(1,1,1).
\]
Their option sets are
\[
\begin{aligned}
f((0,0,0))
  &=\varnothing,\\
f((1,0,0))
  &=\{(0,0,0)\},\\
f((1,1,0))
  &=\{(0,1,0),(1,0,0)\},\\
f((1,1,1))
  &=\{(0,1,1),(1,0,1),(1,1,0)\}.
\end{aligned}
\]
Therefore, the recursive definitions of $g^+$ and $g^-$ give
\[
\begin{array}{c|cccc}
\m
 &(0,0,0)&(1,0,0)&(1,1,0)&(1,1,1)\\ \hline
g^+(\m)&0&1&0&1\\
g^-(\m)&1&0&1&0
\end{array}
\]
and the assertion follows by cyclic symmetry.
\end{proof}

\begin{lemma}
\label{lemma_tame_prepare}
Let $\m\in (M\setminus A)$. Then the following conditions are
equivalent:
\begin{enumerate}
    \item[(\RMi)]
    There exists an option $\m'\in f(\m)$ such that
    $\m'\in S_1^+$.
    
    \item[(\RMii)]
    There exists an option $\m'\in f(\m)$ such that
    $\m'\in S_1^-$.
\end{enumerate}
\end{lemma}
\begin{proof}
We write
$$\m=(x,y,z).$$
By cyclic symmetry, we may assume that
$$x=\max(\{x,y,z\}).$$
Since $\m\notin A$, we have $x\geq2$.

We divide into the following cases.

\begin{enumerate}
    \item[(1)] $y\geq2$ or $z\geq2$;
    \item[(2)] $y=1$ and $z\leq1$;
    \item[(3)] $y=0$ and $z\leq1$.
\end{enumerate}

We show that (\RMi) $\Leftrightarrow$ (\RMii) for these cases one by one.

\begin{enumerate}
    \item[(1)]  Since $x\geq2$, there is no option of $\m$ belongs to $A$. Therefore,
    both $(\RMi)$ and $(\RMii)$ are false.
    \item[(2)] Since $x\geq2$, the following two moves are legal:
    \[
    \m=(x,1,z)\longrightarrow \m'_1=(1,1,z)
    \]
    and
    \[
    \m=(x,1,z)\longrightarrow \m'_0=(0,1,z).
    \]
    Both $\m'_1$ and $\m'_0$ belong to $A$, and
    $
    \ell(\m'_1)=\ell(\m'_0)+1.
$
    Hence, $\ell(\m'_1)$ and $\ell(\m'_0)$ have opposite
    parity. Thus, one of $\m'_1$ and $\m'_0$ belongs to
    $S_1^+$, and the other belongs to $S_1^-$.
    Therefore, both $(\RMi)$ and $(\RMii)$ hold.
    \item[(3)]Since $x\geq2$, the following two moves are legal:
    \[
    \m=(x,0,z)\longrightarrow \m'_1=(1,0,z)
    \]
    and
    \[
    \m=(x,0,z)\longrightarrow \m'_0=(0,0,z).
    \]
    Both $\m'_1$ and $\m'_0$ belong to $A$. 
    As in case~(2), one of $\m'_1$ and $\m'_0$ belongs to
$S_1^+$, and the other belongs to $S_1^-$.
Therefore, both $(\RMi)$ and $(\RMii)$ hold.
\end{enumerate}
\end{proof}

\begin{theorem}
\label{thm:triangle_tame}
For every $\m\in M\setminus A$,
\[
g^+(\m)=g^-(\m).
\]
Consequently, Triangle Game is tame.
\end{theorem}

\begin{proof}
We proceed by induction on $\ell(\m)$, in the form of a minimal-counterexample argument.

Suppose, to the contrary, that there exists a position
$\m\in M\setminus A$ such that
\[
g^+(\m)\neq g^-(\m).
\]
We choose such a position $\m$ for which $\ell(\m)$ is minimum.

Since every move strictly decreases the total
number of stones, every option $\m'\in f(\m)$ satisfies
$
\ell(\m')<\ell(\m).
$
Therefore, by the minimality of $\ell(\m)$,
\[
g^+(\m')=g^-(\m')
\]
for every option $\m'\in f(\m)\setminus A$.

First, suppose that $\m$ has no option in $A$, that is,
\[
f(\m)\cap A=\varnothing.
\]
Hence
\[
\{g^+(\m')\mid \m'\in f(\m)\}
=
\{g^-(\m')\mid \m'\in f(\m)\}.
\]
Since $\m\notin A$, we have $\m\neq(0,0,0)$. Taking the mex of the two equal sets gives
\[
g^+(\m)=g^-(\m),
\]
which is a contradiction.

Next, suppose that $\m$ has an option in $A$, that is,
\[
f(\m)\cap A\neq\varnothing.
\]
Since
$
A=S_1^+\cup S_1^-,
$
the position $\m$ has an option in at
least one of $S_1^+$ and $S_1^-$. By
Lemma~\ref{lemma_tame_prepare}, it therefore has an option in
each of these two sets.

By Lemma~\ref{gyakuten}, every position in $S_1^+$ is a
$(0,1)$-position, whereas every position in $S_1^-$ is a
$(1,0)$-position. Therefore,
\[
\{g^+(\m')\mid \m'\in f(\m)\cap A\}
=
\{0,1\},
\]
and
\[
\{g^-(\m')\mid \m'\in f(\m)\cap A\}
=
\{0,1\}.
\]

Consequently,
\begin{align*}
g^+(\m)
&=
\Mex\bigl(
 \{g^+(\m')\mid \m'\in f(\m)\}
\bigr)\\
&=
\Mex\bigl(
 \{0,1\}
 \cup
 \{g^+(\m')\mid
   \m'\in f(\m)\setminus A\}
\bigr)\\
&=
\Mex\bigl(
 \{0,1\}
 \cup
 \{g^-(\m')\mid
   \m'\in f(\m)\setminus A\}
\bigr)\\
&=
\Mex\bigl(
 \{g^-(\m')\mid \m'\in f(\m)\}
\bigr)\\
&=
g^-(\m).
\end{align*}
Here, the third equality follows from the minimality of
$\ell(\m)$, and the fourth equality follows from
Lemmas~\ref{gyakuten} and \ref{lemma_tame_prepare}.
This again contradicts the choice of $\m$. 

Finally, by Lemma~\ref{gyakuten}, every position in $A$ is
either a $(0,1)$-position or a $(1,0)$-position. Every position
outside $A$ is a $(k,k)$-position for some
$k\in\mathbb Z_{\geq0}$. Therefore, Triangle Game is tame.
\end{proof}

\section{Final Remarks}
\label{Sec:final}

We have shown that the $\mathcal{P}$-positions of Triangle Game can be characterized in terms of the golden ratio. Games whose $\mathcal{P}$-positions admit such a characterization are quite rare, and we believe that this result is significant in revealing a new connection between combinatorial games and the golden ratio.

One possible direction for future research is to generalize the present game. 
Let $G=(V,E)$ be a directed graph, and assign a
non-negative number of stones to each vertex. For a directed edge $U\to V$,
one may consider a move in which a player removes $i\ge1$ stones
from $U$ and transfers $j$ of the removed stones to $V$, subject
to
$$
0\le j<i.
$$
Triangle Game is one of the special cases in which \(G\) is a directed
$3$-cycle.

It is natural to ask for which directed graphs the corresponding
$\mathcal P$-positions admit an explicit description. In preliminary
investigations, we considered this game on several other directed
graphs, including directed cycles of
length greater than three. We also considered variants obtained by modifying the
restriction $0\le j<i$.
See~\cite{IKMSYY} for studies on such a variant, including the computational complexity of its winner-determination problem.

However, we have not yet found another example in which the
$\mathcal P$-positions can be characterized in a way similar to that used for
Triangle Game using the golden ratio or another algebraic number.

This leads to the following questions. Which structural properties of
a directed graph allow a simple description of the
\Pps? Under what conditions does an irrational
number arise naturally in such a description? 

Another direction for future research is to determine not only the $\mathcal{P}$-positions but also the Sprague--Grundy values. 
From our results, it is straightforward to verify that the Triangle Game is tame, meaning that it satisfies a certain relationship between its normal and misère Sprague--Grundy values. Unfortunately, however, we were unable to find a simple closed-form expression for the Sprague--Grundy value of an arbitrary position of the Triangle Game. This suggests that the problem is highly nontrivial.

\section*{Acknowledgments}
The authors would like to thank the anonymous readers for their valuable comments and suggestions.


\end{document}